\documentclass[thmsa,notitlepage,11pt]{article}
\usepackage{sw20bams}

\input tcilatex
\QQQ{Language}{
American English
}

\begin{document}

\title{An Algebraic Treatment of Totally Linear \\
Partial Differential Equations}
\author{C. Viazminsky \\
International Institute of Theoretical and Applied Physics,\\
\ \ \ Iowa State University, Ames, IA 50011, USA \\
and Department of Physics, Aleppo University, Syria}
\maketitle

\begin{abstract}
We construct the field generated by $n$ algebraically independent elements,
and show that the linear space of derivations over this field is faithfully
represented by the linear space of the $n-th$ fold Cartesian product of this
field acting through inner product on the gradient of this field. We prove
also that functional independence of a set in this field is equivalent to
linear independence of the gradient set in the vector space of Cartesian
product. It is also shown that every subfield $S$ of $A$ which is generated
by $(n-1)$ functionally independent elements defines a one- dimensional
space of derivations, such that each member $L$ of the latter subspace has $%
S $ as its kernel. Each coset of the multiplicative subgroup $S$ defines a
non-homogeneous differential operator $L+q$ whose kernel coincides with this
coset. We prove also that every element in $A$ defines a coset of the
subgroup $Ker(L+q)$ in the additive group $A$, on which $L+q$ is constant.
\end{abstract}

\section{Introduction}

In a recent work \cite{Viazminsky} the relations between the solutions of
three types of partial differential equations involving first order
differential operators were obtained. It was shown that if $L$ is a
homogeneous differential operator of class $C^0$ on a manifold $M$, $q$ is $%
C^{0\text{ }}$function on $M$ then all operators of the form $L+q$ are
isomorphic to each other when acting on appropriate Hilbert spaces of
functions on $M$. Neat relations between the solutions of the totally linear
partial differential equations$\;(i)\;L\phi =0,\;(ii)\;(L+q)\psi
=0,\;(iii)\;(L+q)\chi =b,$ where $b$ is a $C^0$ function, were derived. More
precisely, it was shown that if $\eta \in Ker(L+q),$ is any particular
solution of $(ii)$ then the set of solutions of $(ii)$ is given by $%
Ker(L+q)=\eta KerL.$ If $\zeta $ is any solution of $(iii)$ then the set of
solutions of $(iii)$ is given by $\zeta +Ker(L+q).$ The set of solutions of $%
(i)$ is a sub-algebra of $C^1(M)$, whereas the set of solutions of $(ii)$, $%
Ker(L+q)$, is a sub-space of the vector space $C^1(M)$. The algebraic
features of these results encourages us to seek an algebraic approach to
this problem. Following up this line of thoughts

- we first formulate the problem in a field $A$ generated by a finite number
of elements, say $n$ elements.

- study derivations in the algebra defined naturally by this field, and show
that the set of derivations is an $n$-dimensional vector space over the
field $A$. In this concern we come close to the treatment of vector fields
in manifold theory \cite{Bishop}.

- show that functional independence of a set in the field $A$ is equivalent
to linear independence of the gradient of this set in the vector space $A^n$.

- Utilize the latter result to find the maximum number of functionally
independent elements in $KerL.$

- Apply the latter results to the solution of equations $(ii)$ and $(iii)$.

- show that every subfield $S$ of $A$, which is generated by $(n-1)$
functionally independent elements defines a derivation whose kernel is $S,$
and that every coset of this multiplicative subgroup defines a
non-homogeneous differential operator whose kernel is this coset.

\section{Algebraic Setting}

Let $A$ be a field of characteristic $0,$ and let $K$ be the sub-field
generated by the identity element $1\in A.$ The field $K$ is isomorphic to
the field of rational number $Q$, and will be called the scalar field. The
field $A$ is a vector space over the scalar field $K$, and hence it is a
commutative algebra $(A,K)$ with an identity $1$ and every element in $%
(A,K), $ that is different of zero, is invertible.

We assume that the field $A$ is generated by a set $G=\{x_1,...,x_n\}\;$%
consisting of $n$ distinct elements that are different from the field's
zero, $0$, and identity,$1$, and that the set $G$ is functionally
independent in the sense that there is no conceivable polynomial in the
elements of this set equals to zero except the zero polynomial (appendix 1).
It follows from definition of a generating set that every element $\psi \in
A $ is a function of the elements of the generating set $G$. A typical
element $\psi \in A$ is of the form 
\begin{equation}
\psi =P_1P_2^{-1},  \label{e1}
\end{equation}
where $P_1$ and $P_2\neq 0\ $are two arbitrary polynomials in elements of G.
More details about the field $A$, the correlation between algebraic and
functional independence, and how to construct a basis in the infinite
dimensional vector space $(A,K)$ is given in appendix 1.

We assume also that the field $A$ is partially ordered by a partial ordering
relation $(\leq )$ that has the following additional properties

$(i)$\ \ $\psi _1>\psi _2,\;\varphi _1\geq \varphi _2\Rightarrow \psi
_1+\varphi _1>\psi _2+\varphi _2$ .

$(ii)$\ \ $\psi ^2\geq 0.$

From $(ii)$ and the fact that every element in $A-\{0\}$ is invertible, we
get

$(iii)$ \ $\psi ^2=0\Leftrightarrow \psi =0.$

From $(i)$ and $(iii)$, it follows that

$(iv)$\ \ $\psi ^2+\varphi ^2=0\Leftrightarrow \psi =0,\;\varphi =0.$

Consider the Cartesian product $A^n.$ The set $A^n$ forms an $n$-dimensional
vector space over the field $A$, when addition and multiplication by a
scalar from $A$ are defined coordinate wise. We shall denote this vector
space by $(A^n,A)$. The mapping $\prec .\mid .\succ :A^n\times
A^n\rightarrow A,$ defined by 
\begin{equation}
\prec (a_1,....,a_n)\mid (b_1,....,b_n)\succ =\stackunder{i=1}{\stackrel{n}{%
\sum }}a_ib_i  \label{e2}
\end{equation}
is left linear, symmetric, and positive definite, and hence is an inner
product in $A^n.$ The property of positive definiteness can be checked using
the additional property $(iv)$ of the partial ordering relation.

\section{The Space of Derivations}

Let $L:(A,K)\rightarrow (A,A)$ be a linear mapping that satisfy the
condition 
\begin{equation}
L(\psi \varphi )=\psi (L\varphi )+\varphi (L\psi ).  \label{e3}
\end{equation}
In other words $L$ is a derivation \cite{Matsushima}, \cite{Conlon} from the
algebra $(A,K)$ to the algebra $(A,A)$. It follows from definition that 
\begin{equation}
L\alpha =0\;\;\forall \alpha \in K,\;\;\;L\psi ^n=n\psi ^{n-1}L\psi
\;\;\forall n\in Z.  \label{e4}
\end{equation}
The set of derivations from $(A,K)$ to $(A,A)$ forms a vector space $H$ over
the field $A$, which we denote by $(H,A)$, and conveniently call the space
of homogeneous first order differential operators in $(A,K)$. Let $%
B=\{\partial _1,...,\partial _n\}\ $be a subset of $(H,A)$ satisfying the
property 
\begin{equation}
\partial _ix_j=\delta _{ij\;}\;(i,j=1,...,n).  \label{e5}
\end{equation}
We shall prove that the set $B$ is a basis of the vector space $(H,A)$, and
hence the dimension of the latter space is $n$. Our proof is quite similar
to that followed in manifold theory when proving that the set $B$ is a basis
for the tangent space at some point of a manifold \cite{Bishop}. Let $L\in
(H,A)$ and assume that $Lx_i=a_i(i=1,...,n)$. Now, if $\psi \in A$ is
arbitrary then it has the form (\ref{e1}). Acting by $L$ on $\psi ,$ we show
by a straight forward calculation, that 
\begin{equation}
L\psi =a_1\partial _1\psi +.....+a_n\partial _n\psi .  \label{e6}
\end{equation}
Hence 
\begin{equation}
L=a_1\partial _1+....+a_n\partial _n,  \label{e7}
\end{equation}
and $B$ is a generating set in $(H,A)$. To show that $B$ is linearly
independent in the vector space $(H,A)$, we consider the linear combination 
\begin{equation}
\stackrel{n}{\stackunder{i=1}{\sum }}a_i\partial _i=0\;\;\;\;(a_i\in A)
\label{e8}
\end{equation}
Acting by both sides of the latter equality on $x_k$ we get $a_k=0,$ which
is true for all $k\in [1,n].$ Hence the set $B$ is linearly independent, and 
$B\;$is a basis in $(H,A)$.

Note that every $L\in (H,A)\;$is uniquely defined by its values at the
elements of the generating set $G$, and that varying these values in every
possible way leads to the set of all possible derivations from $(A,K)$ to $%
(A,A)$ which are consistent with (\ref{e5}). We note also that one could
have set $\partial _ix_j=\delta _{ij}f_i(x)$ to obtain consequently $L=\sum
a_if_i^{-1}\partial _i.$ Though this not consistent with the familiar
calculus, still it may be interesting to study the unfamiliar derivations
and the associated differential equations. However we shall not consider
this topic in our present work.

We derive one more rule of differentiation, namely that concerned with the
chain differentiation. Let $\phi $ be a function from $A^n$ to $A$, where $%
x=(x_1,...,x_n)\rightarrow \phi (x)$. The set of values of all functions $%
\psi (\phi (x))$ forms a subfield $A_\phi $ of the field $A$ , and a
sub-algebra $(A_\phi ,K)$ of $(A,K)$ . This subfield is generated by a
single element $\phi $. A typical element of $A_\phi $, is of the form $\psi
(\phi )=P_m(\phi )(P_l(\phi ))^{-1},$where $P_m(\phi )=\stackunder{}{%
\stackunder{r=1}{\stackrel{m}{\sum }}\alpha _r}\phi ^r,(\alpha _r\in K).$ It
follows that the space of derivation from $(A_\phi ,K)$ to $(A_\phi ,A_\phi
) $ is one dimensional, and hence all derivations on this space have the
form $f(\phi )D_\phi ,$ where $f(\phi )\in A_\phi ,$ and $D_\phi \phi =1.$

Though a bit tedious, it is straight forward to check that 
\[
\partial _i\psi (\phi (x))=D_\phi \psi (\phi )\partial _i\phi (x) 
\]
The last relation can be checked by substituting for $\phi $ the general
expression (\ref{e1})of an element in $(A,K)$. The last relation shows that 
\[
(\partial _i\psi )(\partial _j\phi )=(\partial _j\psi )(\partial _i\phi
)\;\;\;\;(i,j=1,...,n), 
\]
or 
\[
rank\left[ 
\begin{array}{ccc}
\partial _1\psi & ... & \partial _n\psi \\ 
\partial _1\phi & ... & \partial _n\phi
\end{array}
\right] <2, 
\]
which expresses the fact that $\psi $ and $\phi $ are functionally dependent.

In a similar way the set of functions of the form $\psi =\psi (\phi
_1,...,\phi _r),$ $(r\leq n)$, where $\phi _i(x),\;(i=1,...,r)$ are
functionally independent, forms a subfield $A_{1...r\text{ }}$ of $A$ which
is generated by $\{\phi _i(i=1,...,r)\},$ and a sub-algebra of the algebra $%
(A,K)$. The space of derivations in this sub-algebra is $r-$dimensional,
with $\{D_1,...,D_r\}$, where $D_i\phi _j=\delta _{ij},$ as a basis. It can
be checked that 
\begin{equation}
\partial _i\psi (\phi _1(x),...,\phi _n(x))=(D_1\psi )(\partial _i\phi
_1)+....+(D_r\psi )(\partial _i\phi _r).  \label{e9}
\end{equation}
From the latter relation and (\ref{e7}) it follows that 
\begin{equation}
L\psi (\phi _1(x),...,\phi _r(x))=(D_1\psi )(L\phi _1)+.....+(D_r\psi
)(L\phi _r).  \label{e10}
\end{equation}

\section{Link Between Functional and Linear Independence}

Let $\nabla \equiv (\partial _1,....,\partial _n):(A,K)\rightarrow (A^n,A)$
be defined by 
\begin{equation}
\nabla \phi =(\partial _1\phi ,....,\partial _n\phi )\;\;\;\;\phi \in A.
\label{e11}
\end{equation}
The operator $\nabla $, to be called gradient, is linear and satisfy the
identity $\nabla (\psi \phi )=\phi (\nabla \psi )+\psi (\nabla \phi ).$
Hence $\nabla $ is a derivation from $(A,K)$ to $(A^n,A).$ If $\Phi =\{\phi
_1,..,\phi _r\}\subset A$, where $(r\leq n),$ then the vectors $\nabla \phi
_1,...,\nabla \phi _r\in (A^n,A).$ We shall show that the functional
independence of the set of elements $\Phi $ in the field $A$ is equivalent
to the linear independence of the set vectors $\nabla \Phi =\{\nabla \phi
_1,....,\nabla \phi _r\}$ in the vector space $(A^n,A)$. We shall call $%
\nabla \Phi $ the gradient of the set $\Phi .$ Consider the functional
relation 
\begin{equation}
\psi (\phi _1,....,\phi _r)=0.  \label{e12}
\end{equation}
We have then 
\begin{equation}
\stackrel{r}{\stackunder{j=1}{\sum }}(\partial _i\phi _j)(D_j\psi
)=0\;\;(i=1,...,n).  \label{e13}
\end{equation}
Note first that $D_j\psi =0\in (A^r,A)$ is equivalent to $\psi =\alpha $,
where $\alpha \in K.$ However, by (\ref{e12}), we find that $D_i\psi
=0\Leftrightarrow \psi =0.$ i.e. $D_i\psi =0$ is equivalent to the
functional independence of the set $\Phi .$ Equation (\ref{e13}) asserts
that the $r-$vector $(D_j\psi )\in (A^r,A)$ is orthogonal to the $r-$vectors 
$\partial _i\phi _j$ $(i=1,...,n)$. If the $r-$vector $(D_j\psi )$ is not
zero then it can be orthogonal to $(r-1)$ vectors in $(A^r,A)$ at most.
Hence the rank of the $(n\times r)$ matrix ( $\partial _i\phi _j)$ is less
or equal to $(r-1)$, and the set of vectors $\{\nabla \phi
_j:j=1,...,r\}\subset (A^n,A)$ is linearly dependent. Conversely, and if the
set of vectors $\nabla \phi _j$ $\in (A^n,A),$ $(j=1,...,r)$ is linearly
dependent then the rank of the matrix $(\partial _i\phi _j)$ is less than $%
r. $ As a result equation (\ref{e13}) can have a non-trivial solution $%
(D_j\psi )\neq 0.$ Therefore the set of functions $\Phi $ is functionally
dependent.

It is clear that if $r>n$, then the set of functions $\Phi $ is certainly
functionally dependent. This is because the system of linear equations (\ref
{e13}) is under determined, and has consequently a non-trivial solution $%
(D_1\psi ,....,D_r\psi )$. We state the latter results in the following
theorem

\begin{proposition}
1. A subset $\Phi =$ $\{\phi _1,...,\phi _r\}$ of the field $A$ is
functionally independent if and only if the set of vectors $\nabla \Phi
=\{\nabla \phi _1,....,\nabla \phi _r\}\subset (A^n,A)$ is linearly
independent. In other words, the set $\Phi $ is functionally independent iff
its gradient is linearly independent.
\end{proposition}

2. If $r>n$ then $\Phi $ is functionally dependent. Hence there can exists
in $A$ no more than $n$ functionally independent elements. If $r\leq n$ then 
$\Phi $ is functionally independent if and only if 
\begin{equation}
rank\;(\partial _i\phi _j)=r\;\;\;\;(i\in [1,n],\;j\in [1,r])  \label{e14}
\end{equation}

\section{The Homogeneous Linear Equation $L\phi =0.$}

Let $L=\stackrel{n}{\stackunder{i=1}{\sum }}a_i\partial _i$ $\in (H,A).$ The
correspondence 
\begin{equation}
L\in (H,A)\leftrightarrow (a_1,....,a_n)\in (A^n,A)  \label{15}
\end{equation}
is a vector space isomorphism. Consider a subset $\nabla \;A$ of $A^n,$
defined by 
\begin{equation}
\nabla \;A=\{\nabla \phi :\phi \in A\}.  \label{e16}
\end{equation}
The set $\nabla \;A$ will be called the gradient of the field A. It is clear
that $\nabla \;A$ is a subspace of the vector space $(A^n,K).$ Since 
\[
L\phi =\left\langle (a_1,...,a_n)\mid \nabla \phi \right\rangle 
\]
we may identify $L$ by the operator 
\begin{equation}
\left\langle (a_1,...,a_n)\mid \;\right\rangle :(\nabla \;A,K)\rightarrow A
\label{e17}
\end{equation}

If $\phi _k$ is a solution of the homogeneous differential equation ($HDE$
for short) 
\begin{equation}
L\phi =0,  \label{e18}
\end{equation}
then 
\begin{equation}
\left\langle (a_1,...,a_n)\mid \nabla \phi _k\right\rangle =0.  \label{e19}
\end{equation}
It follows that $L\phi _k=0$ if and only if $\;\nabla \phi _k\in L^{\perp },$
where $L^{\perp }\subset (A^n,A)$ is the set of all vectors in $(A^n,A)$
which are orthogonal to the vector $(a_1,....,a_n)\in (A^n,A).$ The set $%
L^{\perp \text{ }}$ is clearly an $(n-1)-$dimensional subspace of the vector
space $(A^n,A)$. Since $L\phi _k=0\Leftrightarrow \phi _k\in \ker L$, we
have $\nabla \ker L\subset L^{\perp }.$ Thus $\nabla \ker L$ can have at
most $(n-1)$ elements, which are linearly independent, and hence $kerL$ can
have at most $(n-1)$ elements that are functionally independent. Therefore
the $HDE$ (\ref{e18}) can have at most $(n-1)$ functionally independent
solutions, which is a well-known result in partial differential equations.

\section{Non-Homogeneous Operators}

Let $q\in A,$ and consider the operator $L+q:(A,K)\rightarrow (A,A),$
defined by 
\begin{equation}
(L+q)\psi =L\psi +q\psi .  \label{e20}
\end{equation}
This operator is linear. The set of such operators $NH=\{L+q:L\in H\;,q\in
A\}$ form a vector space over the field $A$. Through the inclusion mapping $%
L\in (H,A)\rightarrow (L+0)\in (NH,A),$ where $0\in A,$ we may write $%
(H,A)\subset (NH,A).$ We shall show that the set $NB=\{\partial
_1,...,\partial _n,1\}\subset (NH,A)$ is a basis in $(NH,A),$ and thus $%
(NH,A)$ is an $(n+1)$-dimensional vector space. It clear first that $NB$
generates the vector space $(NH,A)$. Consider next the linear combination 
\[
\stackrel{n}{\stackunder{i=1}{\sum }}a_i\partial _i+b.1=0
\]
Acting by both sides on $\alpha \in K$ we get $b\alpha =0$ for every $\alpha
,$ and hence $b=0$. Acting on $x_k\;(k=1,...,n)$ we get $a_k=0\;(k=1,...,n).$
Therefore the set $NB$ is linearly independent, and hence is a basis of $%
(NH,A)$.

An alternative way of reaching the latter results is based on identifying $%
(NH,A)$ with $(H,A)\times (A,A)$ through the isomorphism $L+q\in
(NH,A)\longleftrightarrow (L,q)\in (H,A)\times (A,A).$

If $L=\sum a_i\partial _i$ then the correspondence 
\begin{equation}
L+q\in (NH,A)\longleftrightarrow (a_1,...,a_n,q)\in (A^{n+1},A)  \label{e21}
\end{equation}
is a vector space isomorphism. Since 
\begin{equation}
(L+q)\psi =\left\langle (a_1,...,a_n,q)\mid (\;\nabla \psi ,\psi
)\right\rangle \in A,  \label{22}
\end{equation}
we may identify $L+q$ by the operator 
\begin{equation}
\left\langle (a,...,a,q)\mid (.,.)\right\rangle :(\;\nabla \;A\times
A)\rightarrow A.  \label{23}
\end{equation}
Now $\psi \in A$ is a solution of the non-homogeneous differential equation (%
$NHDE$ for short) 
\begin{equation}
(L+q)\psi =0,  \label{e24}
\end{equation}
if and only if 
\[
\left\langle (a_1,...,a_n,q)\mid (\nabla \psi ,\psi )\right\rangle =0. 
\]
Therefore the set of solutions of the $NHDE$ (\ref{e24}), namely $Ker(L+q)$
satisfy the condition 
\begin{equation}
\nabla Ker(L+q)\times Ker(L+q)\subset (L+q)^{\perp },  \label{e25}
\end{equation}
where $(L+q)^{\perp }$ is the set of all vectors in $(A^{n+1},A)$ which are
orthogonal to the vector $(a_1,...,a_n,q)\in (A^{n+1},A).$ Since the set $%
(L+q)^{\perp }$ is an $n-$dimensional sub-space of the vector space $%
(A^{n+1},A)$, we conclude that there exist at most $n$ linearly independent
vectors $(\nabla \psi _k,\psi _k)\in \nabla \ker L\times \ker L$ $%
(k=1,...,n) $. Equivalently there exist at most $n$ solutions of the $NHDE$ (%
\ref{e24}) such that 
\begin{equation}
rank(\partial _i\psi _k,\psi _k)=n.  \label{e26}
\end{equation}
When the last equation holds, it is satisfied in particular, by $n$
solutions $\psi _k$ such that $rank(\partial _i\psi _k)=n.$

\section{Link Between the Solutions of $HDEs$ and $NHDEs$}

The set of solutions of the $HDE$ (\ref{e18}), namely $kerL$, is a
sub-algebra of the algebra $(A,K)$. For, $kerL$ is a subspace of the vector
space $(A,K)$ because $L$ is linear, and $kerL$ is closed under
multiplication because $L$ is a derivation. The set of solutions of the $%
NHDE $ (\ref{e24}) is a subspace of the linear space $(A,K)$. In what
follows we shall assume that $\eta \in A$ is any particular solution of the $%
NHDE$ (\ref{e24}),i.e. $\eta \in Ker(L+q).$

\begin{proposition}
The set of solutions of the $NHDE$ (\ref{e24}) is given by 
\begin{equation}
Ker(L+q)=\eta KerL.  \label{e27}
\end{equation}
Proof: Let $\phi \in KerL,$ then 
\[
(L+q)(\eta \phi )=\phi (L+q)\eta +\eta L\phi =0.
\]
Therefore $\eta \phi \in Ker(L+q),$ and $\eta KerL\subseteq Ker(L+q).$ Let $%
\psi \in Ker(L+q).$ By equation (\ref{e24}) which is satisfied by $\psi $
and $\eta ,$ we have $\psi L\eta =\eta L\psi .$ Now 
\[
L(\psi \eta ^{-1})=\eta ^{-2}(\eta L\psi -\psi L\eta )=0.
\]
Hence $\psi \eta ^{-1}\in KerL,$ and $\psi =\eta (\psi \eta ^{-1})\in \eta
KerL,$ and the theorem is proved.
\end{proposition}

\begin{corollary}
In the context of the latter proposition we proved that if $\xi $ and $\eta $
are particular solutions of the $NHDE$ (\ref{e24}) then $\xi \eta ^{-1\text{ 
}}$is a particular solution of the $HDE$ (\ref{e18}). However this result is
contained in equation (\ref{e27}) which gives $\xi \eta ^{-1}KerL=KerL.$ But
since $KerL$ is a subgroup of the multiplicative group $A$, the last
equation implies that $\xi \eta ^{-1}\in KerL.$
\end{corollary}

\begin{corollary}
Equation (\ref{e27}) shows that the general solution of (\ref{e24}) is given
by the product of a particular solution of this equation by the general
solution of equation (\ref{e18}).
\end{corollary}

\begin{proposition}
$\eta \in Ker(L+q)$ $\Leftrightarrow $ $L+q=\eta L\eta ^{-1}$.
\end{proposition}

Proof: It is understood that the latter equation is an operator equation, in
which $\eta L\eta ^{-1}$ is the composite of three mapping in $(A,K)$, so
that $(\eta L\eta ^{-1})\psi =\eta (L(\eta ^{-1}\psi )).$

$(\Rightarrow )$ Assume first that $(L+q)\eta =0,$ so that $q=-\eta
^{-1}(L\eta )$. Hence 
\[
\eta L\eta ^{-1}=L-\eta \eta ^{-2}(L\eta )=L+q. 
\]
$(\Leftarrow )$ $(L+q)\eta =(\eta L\eta ^{-1})\eta =\eta (L1)=0.$ Hence $%
\eta \in Ker(L+q),$ and the proof is complete.$.$

Consider now the $NHDES$ ($S$ is short for (with second side)) 
\begin{equation}
(L+q)\chi =b,  \label{e28}
\end{equation}
where $b\in A.$ If $\chi _0$ is a particular solution of this equation and $%
\chi $ is any other solution, then 
\begin{equation}
(L+q)(\chi -\chi _0)=0,  \label{29}
\end{equation}
and the solution of the $NHDES$ is reduced to the solution of the $NHDE$.
Therefore the set of solution $O$ of (\ref{e28}) is given by 
\begin{equation}
O=\chi _0+Ker(L+q)=\chi _0+\eta KerL.  \label{e30}
\end{equation}
The last relation reads that the set of solutions of a $NHDES$ is a coset of
the subspace $Ker(L+q)$ determined by a particular solution $\chi _0.$

\section{Correspondence between Cosets and Differential Operators}

Let $\Phi =\{\phi _1,...,\phi _{n-1}\}$ be a functionally independent subset
in $A$, and let $S$ be the sub-field generated by $\Phi .$ Since $1\in S,$
the scalar field $K$ is a sub-field of $S$. Evidently, $S$ is a sub-algebra
of $A$. Let $\Psi =\{\psi _1,...,\psi _{n-1}\}$ be any other functionally
independent subset of $S$. Since $\Phi $ generates $S$ we have 
\[
\psi _i(x)=f_i(\phi _1(x),....,\phi _n(x))\;\;\;\;\;(i=1,...,n-1). 
\]
Hence 
\[
\left[ 
\begin{tabular}{l}
$\nabla \psi _1$ \\ 
.... \\ 
$\nabla \psi _{n-1}$%
\end{tabular}
\right] =\left[ 
\begin{tabular}{lll}
$D_1f_1$ & ...... & $D_{n-1}f_1$ \\ 
........ & ...... & ......... \\ 
$D_1f_{n-1}$ & ..... & $D_{n-1}f_{n-1}$%
\end{tabular}
\right] \left[ 
\begin{tabular}{l}
$\nabla \phi _1$ \\ 
..... \\ 
$\nabla \phi _{n-1}$%
\end{tabular}
\right] . 
\]
This relation shows that every vector $\nabla \psi _i(i=1,...,n-1)$ is a
linear combination in the elements of the set $\nabla \Phi .$ Since $\nabla
\Phi $ and $\nabla \Psi $ are both linearly independent in $(A^n,A)$, it
follows that the $(n-1)$-dimensional vector space generated by $\nabla \Phi $%
, denoted by $[\nabla \Phi ]$ is identical to $[\nabla \Psi ]$. We shall
denote this sub-space by $L^{\bot }.$ The orthogonal complement of $L^{\bot
} $ is an one-dimensional space spanned by one vector $(a_1,...,a_n)$. This
determines an one dimensional space of homogeneous differential operators
spanned by $L=\stackrel{n}{\stackunder{i=1}{\sum }}a_i\partial _i$ and such
that $kerL=S.$ It is clear of course that $L$ is any member in this one
dimensional space, and hence $L$ is determined up to a multiplicative
function $f(x)\in A.$ We conclude therefore that every functionally
independent subset of $(n-1)$ elements determines on one hand a sub-field $S$
of $A$, and on the other hand an one dimensional space of homogeneous
differential operators, such that if $L$ is any member of the latter
sub-space then $KerL=S.$

Define now an equivalence relation on $A$ by 
\[
\xi \sim \eta \Leftrightarrow \xi \eta ^{-1}\in S 
\]
This relation partitions $A$ to the subgroup $S$ (subgroup with respect to
multiplication in $A$) and all its cosets. Consider the coset defined by $%
\eta $, where $\eta \notin S.$ Since ($\eta L\eta ^{-1})\eta =0,$ it follows
that $\eta \in \ker (L+q)$, where $q=-\eta ^{-1}(L\eta ).$ If $\xi $ is any
other element in the coset $\eta S$, then $(\eta L\eta ^{-1})\xi =\eta
L(\eta ^{-1}\xi )=0,$ because $\eta ^{-1}\xi \in S=\ker L.$ Therefore every
coset $\eta S$ of $S$ defines, up to a multiplicative element from $A$, a
unique non-homogeneous operator $L+q$ such that $\eta S=\ker (L+q).$

Now every proper coset $\eta \ker L=\ker (L+q)$ of the subgroup $KerL$ is a
subgroup of the additive group $A$. Define on $A$ an equivalence relation by 
\[
\chi \simeq \chi _0\Leftrightarrow \chi -\chi _0\in \ker (L+q). 
\]
This equivalence relation partitions $A$ into equivalence classes consisting
of $Ker(L+q)$ and all its cosets. Take $\chi _0\in A,$ such that $\chi
_0\notin \ker (L+q),$ and set $(L+q)\chi _0=b.$ Let $\chi \in \chi _0+\ker
(L+q)$ be arbitrary. Since $\chi -\chi _0\in \ker (L+q),$ we have $%
(L+q)(\chi -\chi _0)=0,$ and hence $(L+q)\chi =(L+q)\chi _0=b.$ Therefore,
given any operator $(L+q)$, every $\chi _0\in A,$ determines a coset of $%
ker(L+q)$, in the additive group $A$, on which the operator $L+q$ is
constant.

\section{Appendix 1.}

Here we list the ordered set $B_P$ consisting of all one term polynomials
(or the non fractional elements) in $A$ according to raising powers

$1;$

$x_1,..........,x_n;$

$x_1^2,x_1x_{2,}...,x_1x_n;\;x_2^2,x_2x_3,...,x_2x_n;%
\;x_3^2,x_3x_4,....,x_3x_n;....;x_{n-1}^2,x_{n-1}x_{n;}\;x_n^2;$

$x_1^3,x_1^2x_2,...;x_2^3,x_2^2x_3,.....$

$................................$

Now every polynomial $P\in A$ can be written as a linear combination in the
elements of $B_P.$ The assumption that the generating set is functionally
independent is equivalent to say that the only vanishing linear combination
of the elements of $B_P$ is the trivial linear combination. From this
property it follows that every polynomial in $A$ can be expressed in a
unique way as a linear combination of the elements of $B_P.$ Thus to every
polynomial $P\in A$ there corresponds a unique sequence $(s_i)$ in $Q,$
defined by $s_i$ $=$ the $i-th$ coefficient of $P$ on $B_P.$ It clear that
the set of all polynomial in $A$, denoted by $P(A)$ forms an infinite
dimensional vector space over the field $K$, with $B_P,$ as a basis. The
vector space $P(A)$ is isomorphic to the set $S(Q)$ consisting of all
infinite sequences of rational numbers which have only a finite number of
non-zero terms$.$ It is understood of course that addition in $S(Q)$ and
multiplication by a scalar from $Q$ is defined term-wise. As an example,
there corresponds to the polynomial $2^{-1}+y+3.4^{-1}x^2+y^2+x^2y,$ when $%
n=2$, the sequence $(1/2,0,1,3/4,0,1,0,1,0,0.....).$

We mention that the general form of a polynomial of degree $k$ in the
elements of $G$ is 
\[
P=\stackunder{s=0}{\stackrel{k}{\sum }}\stackunder{\;i_1...i_n}{\sum }\alpha
_{i_1.....i_n}^{j_{1........}j_n}\;x_{i_1}^{j_1}....x_{i_n}^{j_n},
\]
where the indices are subject to the conditions 
\[
i_1<...<i_n\leq n,\;\;\;j_1\geq ...\geq j_n\geq 0,\;\;s=j_1+...+j_n\in [0,k],
\]
and where the $\alpha ^{\prime }s$ are scalars from $K.$ Note that any
function of elements of $A$ is a element of $A$ and hence will have the form
(\ref{e1}). If this function is zero then the numerator is zero, and
therefore functional independence in $A$ is equivalent to algebraic
independence. The set $P(A)$ is an integral domain with respect to the usual
operations of addition and multiplication, and it is clear that the field $A$
is the quotient field of the integral domain $P(A)$. The set of sequences $%
S(Q)$ forms also an integral domain, with respect to the usual addition and
multiplication of infinite sequences, which is isomorphic to the integral
domain $P(A)$. Therefore the quotient field of $S(Q)$ is isomorphic to the
quotient field of $P(A)$, which is $A$.

In order to establish a basis for the vector space $(A,K)$, we agree on the
following: we represent every element $\psi =P_1P_2^{-1}$ by a pair $%
(P_1,P_2)\in P(A)\times P(A)^{*},$ such that the coefficient of the leading
term in $P_2$ is equal to $+1$, and such that there is no common factor in $%
P_1$ and $P_2.$ Under this convention one can see that $(B_P,B_P)$ is a
basis of $(A,K)$. It is evident that the addition and multiplication in $A$
are defined by $(P_1,P_2)+(p_1,p_2)=(p_2P_1+p_1P_2,p_1p_2),%
\;(P_1,P_2).(p_1,p_2)=(P_1p_1,P_2p_2)$. If we adopt in $S(Q)\times S(Q)^{*}$
a similar agreement to that adopted in $A$, we find that $A$ is isomorphic
to $S(Q)\times S(Q)^{*}.$

An alternative to the content of this appendix is to consider each element $%
\psi =P_1P_2^{-1}$ as a class of equivalence in $P(A)\times P(A)^{*}$, in
which $(P_1,P_2)\sim (P_3,P_4)$ iff $P_1P_3=P_2P_4,$ and to define a similar
equivalence relation on $S(Q)\times S(Q)^{*},$ so that the quotient field of 
$S(Q)$ is isomorphic to $A$.

The algebraic notions used in this appendix can be found in \cite{Dummit}

\section{Gratitude}

The author thanks Professors Irvin Hentzel and Alex Abian from the
mathematics department in Iowa State University for useful discussions, and
thanks the university of Aleppo in Syria for financial support.

\end{document}